\documentclass[12pt,a4paper,twoside]{article}

\usepackage[T1]{fontenc}
\usepackage[cp1250]{inputenc}
\usepackage[polish,english]{babel}
\usepackage{graphicx}
\usepackage[intlimits]{amsmath}
\usepackage{amsthm}
\usepackage{amsfonts}

\newcommand{\dow}{{\noindent\bf Proof.}}

\newtheorem{lemat}{Lemma}
\newtheorem{twr}[lemat]{Theorem}

\newtheorem{remark}[lemat]{Remark}

\title{Analysis of nonlocal model of compressible fluid in 1-D}
\author{Ewelina Kamińska}
\begin{document}
\maketitle
\normalsize

\noindent{\bf Abstract:} The compressible barotropic Navier-Stokes type system in monodimensional case with Neumann boundary condition given on free boundary is considered. The local and the global existence with uniformly boundedness for small viscosity coefficient is proved.\\
\noindent{\bf Keywords:} {the Navier-Stokes equations, barotropic compressible viscous fluid,
 weak solution, global existence}

\vspace{0.2cm} \noindent{\it MOS subject
classification:} 76N10; \vspace{0.2cm} 35Q30
 \renewcommand{\theequation}{1.\arabic{equation}}
 \setcounter{equation}{0}
 \section{Introduction}
\noindent In this article we consider a model of motion of isolated volume of a barotropic viscous compressible fluid in monodimensional case with a free boundary given by an initial-value problem for modified Navier-Stokes system. This system can be treated as a model of a single layer of a star. The equations in Euler's co-ordinates are of this form, because we want operator $\mathbf{T}$ to act only on the velocity function after transformation to Lagrangian mass co-ordinates, thus
\begin{equation}\label{euler}
\begin{array}{c}
v_{t}+vv_{r}+\frac{1}{\varrho}p_{r}=\mu\mathbf{T}\left(\frac{1}{\varrho}\left(\frac{v_{r}}{\varrho}\right)_{r}\right)\\
\varrho_{t}+(\varrho v)_{r} = 0\\
\mu\mathbf{T}\left(\frac{v_{r}}{\varrho}\right)-p=-P\quad   \mathrm{for}\quad r=0\quad \mathrm{and}\quad  r=S(t)\\
v(r,0)=v_{0}(r),\quad \varrho(r,0)=\varrho_{0}(r),
\end{array}
\end{equation}
where $v,\  \varrho,\ \mu$ and $P$ are the velocity, the density of the fluid, the positive constant viscosity coefficient and the external constant pressure, respectively; $S(t)$ describes the free boundary, we assume that $S(0)=1$. Function $p=p(\varrho)$ describes the pressure of the fluid as a function of the density.\\
Under the physical constraints function $p(\varrho)$ must satisfy
\begin{equation*}
p(0)=0,\quad p(s_{1})<p(s_{2})\quad  \mathrm{ if}\quad s_{1}<s_{2}. 
\end{equation*}
Moreover, function $G(\cdotp)$ given by the relation
\begin{equation}
p(s)=G'(s)s^{2},\label{defG}
\end{equation}
fulfills a condition
\begin{equation*}
G(s)\geq as^{\gamma-1}
\end{equation*}
for $\gamma>1$ and $a>0$.\\
The classical example of such equation is $p(\varrho)=a\varrho^{\gamma}$ that holds for isentropic processes, but  the use of this model for viscous gas is justified if we assume that the viscosity coefficient is small.\\

The studied system is examined in the Lagrangian mass co-ordinates given by
\begin{displaymath}
x=\int_{0}^{r}{\varrho(r',t)dr'},
\end{displaymath}
and its inverse transformation
\begin{displaymath}
r=\int_{0}^{x}{\xi(y,t)dy},
\end{displaymath}
where $\xi(x,t)=\varrho^{-1}(x,t)$.\\
After this transformation problem (\ref{euler}) reads
\begin{equation}\label{eq:system}
\begin{array}{c}
v_{t}+p(\xi^{-1})_{x}=\mu(\mathbf{T}v)_{xx}\\
\xi_{t}-v_{x}=0 \\
\mu(\mathbf{T} v)_{x}-p(\xi^{-1})\Big|_{x=0}=\mu(\mathbf{T} v)_{x}-p(\xi^{-1})\Big|_{x=1}-P\\
v(x,0)=v_{0}(x),\quad\xi(x,0)=\xi_{0}(x).
\end{array}
\end{equation}
Additionally we make the following assumptions:\\
A1. The external pressure $P>0$.\\
A2. The initial values satisfy
\begin{equation*}
\int_{0}^{1}{v_{0}(x)dx}=0
\end{equation*}
and
\begin{eqnarray*}
\xi_{0}(x)>0,\quad\xi_{0}(0)=\xi_{0}(1),\\
\int_{0}^{1}\varrho_{0}(x)dx=\int_{0}^{1}\xi_{0}^{-1}(x)dx=1,
\end{eqnarray*}
the last condition means that the total mass of the fluid is equal to 1.\\
On the right hand side of the first equation of system (\ref{eq:system}) we have a pseudo-differential operator acting on the velocity function $v$, being a modification of the standard Laplacian. Its definition is based on the properties of space of weak solutions to (\ref{eq:system}) which is the Neumann-boundary problem.
Therefore we immerse the space of weak solutions in $L_{2}(0,1)$ which is considered as the closure of linear combinations of the smooth functions that form a standard base for the Neumann-boundary problem
\begin{equation*}
w_{k}(x)=\frac{\cos(\pi kx)}{\|\cos(\pi kx)\|_{L^{2}(0,1)}}\qquad k=0,1,\ldots, 
\end{equation*}
then we may describe any function $f\in L_{2}(0,1)$ as follows
\begin{displaymath}
f(x)=\sum_{k=0}^{\infty}{f_{k}w_{k}(x)}.
\end{displaymath}
Staying within above notation we define an operator
\begin{displaymath}
\mathbf{T}:L_{2}(0,1)\to L_{2}(0,1)
\end{displaymath}
such that
\begin{displaymath}
\mathbf{T}f(x)=\sum_{k=R+1}^{\infty}{f_{k}w_{k}(x)}.
\end{displaymath}
Operator $\mathbf{T}$ is a projector which omits first $R$ mods of the function. This feature causes that the r.h.s. of the first equation of (\ref{eq:system}) describes the dissipation of the energy only for high fluctuations and does not involve low mods. If the system exhibits only low mods the equations have features of the Euler's system for compressible, inviscid flow; for mods grater than $R$ we have Navier-Stokes equations in one dimension and the dissipation of the energy is proportional to viscosity coefficient $\mu>0$. For this case it has been proved (see \cite{M1}) that the global solutions exist, and that any solution tends to the stationary solution.\\
The objective of this paper is to show a global in time existence of regular solutions to the problem (\ref{eq:system}). The main difficulty is to show the uniformly boudedness of the density $\varrho(x,t)\geq\xi_{-}$. The idea comes from P.B. Mucha and requires an assumption of smallness of viscosity coefficient $\mu$, which is the most interestiong case from the physical point of view. For the sake of Neumann-boundary condition we have a global existance without assuptions of smallness of initial data. In case of Dirichlet-boundary condition smallness of data is necessary however it depends only on $\gamma$ \cite{Y}.\\
But to obtain a global in time existence we need first a local in time existence and then several informations about slutions uniformly in time.
There are some results about local in time existence in a general three-space dimensional case for Navier-Stokes equations with Neumann-boundary condition given on a free boundary \cite{Z}. But in this paper we apply the technique similar to the one from Reference \cite{K}, after noticing that the first equation of the system (\ref{eq:system}) may be stated as follows
\begin{equation}
v_{t}+p(\xi^{-1})_{x}=\mu v_{xx}-\mu((1-\mathbf{T})v)_{xx}\label{laplasjan},
\end{equation}
where
\begin{displaymath}
\mu(1-\mathbf{T})v_{xx}
\end{displaymath}
is an analytic function and norm of it is controlled by the energy bound, thus it may be treated as given one as an external force $f$.\\
\\
To avoid questions about the well posedness of considered problem  in the classical sense we will introduce its weak formulation.\\

\noindent {\bf{Definition (Weak solutions).}}{\it{
We say the pair of functions
\begin{displaymath}
v\in W_{2}^{2,1}((0,1)\times(0,T))\quad and\quad \xi \in L_{\infty}(0,T;H^{1}(0,1))
\end{displaymath}
is a weak solution of the problem (\ref{eq:system}) provided:\\
1. equalities
\begin{equation*}
\begin{array}{c}
(v_{t},\varphi)-(p(\xi^{-1})-P,\varphi_{x})+\mu(\mathbf{T}v_{x},\varphi_{x})=0\\
\xi_{t}-v_{x}=0
\end{array}
\end{equation*}
are fulfilled in the sense of distributions on time interwal [0,T] for each \mbox{ $\varphi \in C^{1}(0,T;H^{1}(0,1))$} of the structure $\varphi(x,t)=\frac{-\pi_{t}(t)}{2}+x\pi_{t}(t)+\eta(x,t)$, where $\int_{0}^{1}{\eta(x,t)dx}=0$ and\\
2. $v(x,0)=v_{0}$, $\xi(x,0)=\xi_{0}$.}}\\
\noindent In above definition we require from the function $v$ regularity, which is not optimal to the weak formulation, however we will need such high smoothness to show uniqueness of the solution.\\
\\
\noindent The results are the following.\\
\\
{\bf{Theorem A (Local in time existence).}}
{\it{Let $v_{0}\in H^{1}(0,1)$, $\int_{(0,1)}{v_{0}(x)dx}=0$ and let $\xi_{0}\in H^{1}(0,1)$, $\xi_{0}>0$. Then there exists $T_{0}>0$ such that there exists unique solution to the problem (\ref{eq:system}) in the sense of Definition on time interval $[0,T]$, $T<T_{0}$.
Moreover, $\xi$ is strictly positive on time interval $[0,T]$.}}\\
\\
\noindent {\bf{Theorem B (Global in time existence).}}
{\it{Let $v_{0}\in H^{1}(0,1)$, $\int_{(0,1)}{v_{0}(x)dx}=0$ and let $\xi_{0}\in H^{1}(0,1)$, $\xi_{0}>0$. Then there exists global in time solution to problem (\ref{eq:system}) such that
\begin{equation}\nonumber
v\in W^{2,1}_{2(loc)}((0,1)\times(0,\infty)),\qquad \xi\in L_{\infty}(0,\infty;H^{1}(0,1)),
\end{equation}
and
\begin{equation}\nonumber
0<\xi_{-}\leq\xi(x,t)\leq\xi_{+}<\infty,
\end{equation}
for all $(x,t)\in(0,1)\times[0,\infty)$, where $\xi_{-},\ \xi_{+}$  are strictly positive constants.}}\\
Notations:
\begin{equation*}
\frac{\partial f(x,t)}{\partial t}=f_{t}(x,t),\qquad \frac{\partial f(x,t)}{\partial x}=f_{x}(x,t),
\end{equation*}
\begin{displaymath}
(f,g)=\int_{0}^{1}{f\cdotp g dx}.
\end{displaymath}

\section{An estimate of the solution on the boundary}\label{rozdz.3}
Note, that from the first equation of the system (\ref{eq:system}) and from the boundary and initial conditions we get
\begin{equation*}
\int_{0}^{1}{v(x,t)dx}=0.
\end{equation*}
Thus it is resonable to look for a function $v(x,t)$ of the structure
\begin{equation}
v(x,t)=\sum_{k=1}^{\infty}\lambda_{k}(t)\omega_{k}(x),\label{szereg}
\end{equation}
where the functions $\lambda_{k}(t)$ are at least $C^{1}(0,\infty)$.\\
This observation together with $\xi_{t}=v_{x}$ let us similarly describe $\xi(x,t)$ as follows
\begin{equation*}
\xi(x,t)=\pi(t)+\sum_{k=1}^{\infty}\kappa_{k}(t)\omega_{k,x}(x),
\end{equation*}
for some function $\pi(t)$ independent of $x$, and the functions $\kappa_{k}(t)$ at least of class $C^{1}(0,\infty)$.
\begin{lemat}\label{nabrzegu}
For all $ t\geq 0$ the function $\xi(x,t)$ satisfies
\begin{equation*}
\begin{array}{c}
\xi(0,t)=\xi(1,t),\\
0<\bar{\xi}_{-}\leq\xi(0,t)\leq\bar{\xi}_{+}<\infty
\end{array}
\end{equation*}
where $\bar{\xi}_{-}$ and $\bar{\xi}_{+}$ are independent of $t$.
\end{lemat}
\dow\\
We examine $\xi(x,t)$ on $(\{0\}\times[0,\infty))\cup (\{1\}\times[0,\infty))$ from the boundary condition of the system (\ref{eq:system})
\begin{equation*}
\mu v_{x} - \mu((1-\mathbf{T})v)_{x}=\frac{a}{\xi^{\gamma}}-P.
\end{equation*}
According to the definition of the operator $\mathbf{T}$ and to the structure of function $v$ we get that for $x=\{0,1\}$ $((1-\mathbf{T})v)_{x}=0$.\\
Therefore the above equation becomes an ODE
\begin{equation}
\mu\xi_{t}=\frac{a}{\xi^{\gamma}}-P\label{brzeg}
\end{equation}
subject to the initial condition  $\xi(x,0)=\xi_{0}(x)$. Thus, locally there exists a unique solution with continous first derivative and it can be extended to the whole half line.\\
From A2 we have $\xi_{0}(0)=\xi_{0}(1)$, hence
\begin{equation*}
\xi(0,t)=\xi(1,t),\quad t\in[0,\infty).
\end{equation*}
Let $\bar{\xi}(t)$ denotes the solution of (\ref{brzeg}) with the initial condition $\bar{\xi}(0)=\bar{\xi}_{0}=\xi(0,0)$.\\
Note, that for $t\rightarrow\infty$
\begin{equation*}
\lim_{t\rightarrow\infty}\bar{\xi}(t)=\left(\frac{a}{P}\right)^{\frac{1}{\gamma}}.
\end{equation*}
Then there are two possibilities:\\
or $\bar{\xi}_{0}\leq\left(\frac{a}{P}\right)^{\frac{1}{\gamma}}$, and the solution is bounded by
\begin{equation*}
0<\bar{\xi}_{-}=\bar{\xi}_{0}\leq\bar{\xi(t)}\leq\left(\frac{a}{P}\right)^{\frac{1}{\gamma}}=\bar{\xi}_{+}<\infty,
\end{equation*}
or $\left(\frac{a}{P}\right)^{\frac{1}{\gamma}}\leq\bar{\xi}_{0}$, and we have 
\begin{equation*}
0<\bar{\xi}_{-}=\left(\frac{a}{P}\right)^{\frac{1}{\gamma}}\leq\bar{\xi(t)}\leq\bar{\xi}_{0}=\bar{\xi}_{+}<\infty. 
\end{equation*}
\begin{flushright}
$\Box$
\end{flushright}
Define the extension $\pi(x,t)$ of the function $\bar{\xi}(t)$ to the whole region $(0,1)\times[0,\infty)$ by the following formula:
\begin{equation*}
\pi(x,t)= \bar{\xi}(t),
\end{equation*}
it means, that te function $\pi(x,t)=\pi(t)$ is constant along each straight line $x=C$.\\

\section{Proof of Theorem A}
First we will construct solutions to some finite-dimensional approximations to (\ref{eq:system}), and then we will pass to the limits. This is called Galerkin's method.\\
\subsection{Galerkin approximations}
Define the spaces:
\begin{multline}\nonumber
W_{v}=\{f\in C^{1}(0,T;H^{1}(0,1)):f=\frac{-\pi_{t}(t)}{2}+x\pi_{t}(t)+h(x,t), \int_{0}^{1}{h(x)dx}=0\}= \\
=\{f\in C^{1}(0,T;H^{1}(0,1)):f=\frac{-\pi_{t}(t)}{2}+x\pi_{t}(t)+h(x,t),\  h=\sum_{k=1}^{\infty}\alpha_{k}(t)w_{k}\},
\end{multline}
and
\begin{equation}\nonumber
W^{N}=\{f\in C^{1}(0,T;H^{1}(0,1)):f=\frac{-\pi_{t}(t)}{2}+x\pi_{t}(t)+h(x,t),\ h=\sum_{k=1}^{N}\alpha^{N}_{k}(t)w_{k}\}.
\end{equation}
For fixed integer $N$, we will look for the functions $v_{N},\ \xi_{N}$ of the form:
\begin{equation}
v_{N}(x,t)=\frac{-\pi_{t}(t)}{2}+x\pi_{t}(t)+\sum_{k=1}^{N}{\alpha^{N}_{k}(t)w_{k}(x)}
\quad\xi_{N}(x,t)=\pi(t)+\sum_{k=1}^{N}{\beta^{N}_{k}(t)w_{k,x}},\label{postac}
\end{equation}
such that for all  $k=1,\ldots N$ the coefficients $\alpha^{N}_{k}(t), \beta^{N}_{k}(t)$ satisfy
\begin{equation}
\alpha^{N}_{k}(0)=(v_{0}-(x-\frac{1}{2})\pi_{t}(0),w_{k}),\quad \beta^{N}_{k}(0)=(\xi_{0}-\pi(0),w_{k,x}),\label{alpha}
\end{equation}
\begin{equation}
\beta^{N}_{k}(t)=\beta^{N}_{k}(0)+\int_{0}^{t}\alpha^{N}_{k}(s)ds,\label{betta}
\end{equation}
and
\begin{equation}\label{gamma}
\begin{array}{c}
(v_{N,t},w_{k})-(\frac{a}{\xi_{N}^{\gamma}}-P,w_{k,x})+\mu(\mathbf{T}v_{N,x},w_{k,x})=0
\end{array}
\end{equation}
in the sense of distributions on time interval $[0,T]$.
\begin{twr}[Construction of approximate solutions]\label{galerkin}
For each integer $N=1,2,\ldots$ there exists the unique pair of functions $v_{N},\ \xi_{N}$ of the form (\ref{postac}) satisfying (\ref{alpha}), (\ref{betta}) and (\ref{gamma}) in the sense of distributions on time interval $[0,T]$.
\end{twr} 
\noindent The Proof is an application of Banach's Fixed Point Theorem.\\
\\
So, the functions $v_{N},\ \xi_{N}$ fulfill the weak formulation in the sense of the Definition for each $\varphi\in W^{N}$.\\
\\
\begin{remark} 
The assumption $\xi_{0}>0$ is equivalent to the initial density $\varrho_{0}(x)=\xi_{0}^{-1}$ bigger than $0$ for all $x\in(0,1)$.\\
This observation, together with the second equation of system (\ref{eq:system}) guarantees 
\begin{equation*}
\xi(x,t)>0\quad \forall(x,t)\in(0,1)\times(0,T),
\end{equation*}
\end{remark}
\noindent Indeed, as
\begin{equation}\nonumber
\xi_{t}(x,t)-v_{x}(x,t)=\left(\frac{1}{\varrho(x,t)}\right)_{t}-v_{x}(x,t)=0,
\end{equation}
hence
\begin{equation}\nonumber
\frac{\varrho_{t}}{\varrho}=-\varrho v_{x},
\end{equation}
thus
\begin{equation}\nonumber
\frac{1}{\xi(x,t)}=\varrho(x,t)=\varrho_{0}(x)\exp\{-\int_{0}^{t}{\varrho v_{x}dt}\}>0.
\end{equation}
This property may be transcribed into $\xi_{N}(x,t)$ (for $N$ sufficiently large) by an analogous argument, we will deal with proving the strict positivity of $\xi_{N}(x,t)$ later.\\
To obtain local existence of weak solutions we will need some uniform estimates.\\
\\
\subsection{Energy estimates}\label{energetyczne}\label{Ienerget}
\begin{lemat}[The first energy estimate]
Let $v_{N},\ \xi_{N}$ satisfy the weak formulation in the sense of the Definition for each $\varphi\in W^{N}$, then
\begin{multline}\nonumber
\int_{0}^{1}{\left(\frac{1}{2}v_{N}^{2}(x,T)+G(\xi_{N}^{-1})(x,T)\right)dx}+PV(T)+\mu\int^{T}_{0}\!\!\!{\int_{0}^{1}{|\mathbf{T}v_{N,x}|^{2}dx}dt}=\\
=\int_{0}^{1}{\left(\frac{1}{2}v_{N,0}^{2}+ G(\xi_{N,0}^{-1})\right)dx}+PR(0)
\end{multline}
holds for any $T< T_{0}$.\\
\end{lemat}
\noindent Putting $v_{N}\in W^{N}$ into the Definition as a test function  $\varphi\in W^{N}$ we get
\begin{equation}
(v_{N,t},v_{N})-(\frac{a}{\xi_{N}^{\gamma}}-P,v_{N,x})+\mu(\mathbf{T}v_{N,x},v_{N,x})=0.\label{zapis1}
\end{equation}
Since $\mathbf{T}$  is the projector, we have
\begin{equation}\nonumber
\mu\int_{0}^{1}{\mathbf{T}v_{N,x}v_{N,x}dx}=\int_{0}^{1}{|\mathbf{T}v_{N,x}|^{2}dx},
\end{equation}
also by (\ref{defG}) 
\begin{equation}\nonumber
-\int_{0}^{1}{p(\xi_{N}^{-1})v_{N,x}dx}=-\int_{(0,1)}{p(\xi_{N}^{-1})\xi_{N,t}dx}=\frac{d}{dt}\int_{0}^{1}{G(\xi_{N}^{-1})dx}
\end{equation}
and
\begin{equation}\nonumber
P\int_{(0,1)}{\xi_{N,t}dx}=P\frac{d}{dt}\int_{(0,1)}{\xi dx}=\frac{d}{dt}PV(t),
\end{equation}
where $V(t)$ is the volume of the fluid in the Eulerian co-ordinates.\\
After this transformations we obtain an equation
\begin{equation}\label{zapis11}
\frac{d}{dt}\int_{(0,1)}{(\frac{1}{2}v_{N}^{2}+G(\xi_{N}^{-1}))dx}+\frac{d}{dt}PV(t)+\mu\int_{0}^{1}{|\mathbf{T}v_{N,x}|^{2}dx}=0,
\end{equation}
integrating from $0$ to $T$ we get the thesis. 
{\begin{flushright}
$\Box$
\end{flushright}}
\noindent Now rewrite equation (\ref{zapis11}) using the formula (\ref{laplasjan})
\begin{equation}
\frac{d}{dt}\int_{0}^{1}{(\frac{1}{2}v_{N}^{2}+G(\xi_{N}^{-1}))dx}+\frac{d}{dt}PV(t)+\mu\int_{0}^{1}{(v_{N,x})^{2}dx}=-\int_{0}^{1}{f_{N}v_{N}dx},\label{pierwsze}
\end{equation}
according to denotation $f_{N}=\mu(1-\mathbf{T})v_{N,xx}=\mu\sum_{k=1}^{R}{\alpha_{k}^{N}(t)w_{k,xx}(x)}$.\\ 
\begin{remark}
Above lemma yields an estimate on the norm of $v_{N}$ in $L_{\infty}(0,T;L_{2}(0,1))$, which in turn implies a suitable estimate on $\|f_{N}\|_{L_{\infty}(0,T;L_{2}(0,1))}$.
\end{remark}
\begin{lemat}[The second energy estimate]\label{IIenerget}
There exists a constant C, depending on $a,\ P,\ \mu,\ \gamma$ and $T$, such that for any $T<T_{0}$ holds
\begin{equation*}
\sup_{t\in[0,T]}\|\xi_{N}(t)\|_{H^{1}(0,1)}+\|v_{N}\|_{L_{2}(0,T;H^{1}(0,1))}\leq C.
\end{equation*}
\end{lemat} 
\dow\\
\noindent This time we take as a test function $\frac{-\pi_{t}(t)}{2}+x\pi_{t}(t)+\xi_{N,x}(x,t)\in W^{N}$, and because  $\frac{-\pi_{t}(t)}{2}+x\pi_{t}(t)$ also belongs to $W^{N}$
we only have to consider
\begin{equation*}
(v_{N,t},\xi_{N,x})-(\frac{a}{\xi_{N}^{\gamma}}-P,\xi_{N,xx})+\mu(\mathbf{T}v_{N,x},\xi_{N,xx})=0.
\end{equation*}
Note
\begin{equation*}
\int_{0}^{1}{v_{N,t}\xi_{N,x}dx}=\frac{d}{dt}\int_{0}^{1}{v_{N}\xi_{N,x}dx}-\int_{0}^{1}{v_{N}\xi_{N,tx}dx}=
\frac{d}{dt}\int_{0}^{1}{v_{N}\xi_{N,x}dx}+\int_{0}^{1}{v_{N,x}^{2}dx}-\pi_{t}\int_{0}^{1}v_{N,x}dx,
\end{equation*}
and
\begin{displaymath}
-\int_{0}^{1}{(\frac{a}{\xi_{N}^{\gamma}}-P)\xi_{N,xx}dx}+\mu\int_{0}^{1}{\mathbf{T}v_{N,x}\xi_{N,xx}}dx=~~~~~~~~~~~~~~~~~~~~~~~~~~~~~~~~~~~~~~~~~~~~~~~~~~~~~~~~~~~~~~~~~~~~
\end{displaymath}
\vskip-18pt
\begin{eqnarray*}
~~~~~~~~~~~~~~~~~~~~~~~~~~~~~~~~~~&=&\int_{0}^{1}{\left(\frac{a}{\xi_{N}^{\gamma}}\right)_{x}\xi_{N,x}dx}-\mu\int_{0}^{1}\xi_{N,xt}\xi_{N,x}dx+\mu\int_{0}^{1}(1-\mathbf{T})v_{N,xx}\xi_{N,x}dx\\
&=&-\gamma\int_{0}^{1}{\frac{a}{\xi_{N}^{\gamma+1}}\xi_{N,x}^{2}dx}-\frac{\mu}{2}\frac{d}{dt}\int_{0}^{1}{\xi_{N,x}^2dx}+\int_{0}^{1}f_{N}\xi_{N,x}dx,
\end{eqnarray*}
thus
\begin{multline}\label{drugie}
\frac{d}{dt}\int_{0}^{1}{\left(\frac{\mu}{2}\xi_{N,x}^2-v_{N}\xi_{N,x}\right)dx}+\gamma a\int_{0}^{1}{\frac{\xi_{N,x}^{2}}{\xi^{\gamma+1}_{N}}dx}=\\
=\int_{0}^{1}{v_{N,x}^{2}dx}+\int_{0}^{1}{f_{N}\xi_{N,x}dx}-\pi_{t}\int_{0}^{1}v_{N,x}dx.
\end{multline}
\noindent Multiply equation (\ref{pierwsze}) by the constant $B=\frac{4}{\mu}$ , and then add to equation (\ref{drugie}), to find
\begin{displaymath}
\frac{d}{dt}\left(\int_{0}^{1}{\left(\frac{\mu}{2}\xi_{N,x}^{2}+\frac{B}{2}v_{N}^{2}-v_{N}\xi_{N,x}+BG(\xi_{N}^{-1})\right)dx}+
BPV(t)\right)+~~~~~~~~~~~~~~~~~~~~~~~~~~~~~~~~~~~~~~
\end{displaymath}
\vskip-18pt
\begin{eqnarray}
~~~~~~~~~~~~~~~~~~~~~~~~~~~~~~~~~~~~~~~~~~&{}&+\gamma a\int_{0}^{1}{\frac{\xi_{N,x}^{2}}{\xi^{\gamma+1}_{N}}dx}+(B\mu-1)\int_{0}^{1}{v_{N,x}^{2}dx}\nonumber\\
&=&\int_{0}^{1}{(f_{N}\xi_{N,x}-Bf_{N}v_{N})dx}-\pi_{t}\int_{0}^{1}v_{N,x}dx.\label{oszacowanie}
\end{eqnarray}
\noindent Employing to the right hand side of  above equation H\"{o}lder's and Cauchy's inequalities and taking into account that $\frac{4}{\mu}\int_{0}^{1}G(\xi_{N}^{-1})>0$, $\frac{4}{\mu}PV(t)>0$ we get
\begin{multline}
\frac{d}{dt}\left(\int_{0}^{1}{\left(\frac{\mu}{2}\xi_{N,x}^{2}+\frac{2}{\mu}{\mu}v_{N}^{2}-v_{N}\xi_{N,x}+\frac{4}{\mu}G(\xi_{N}^{-1})\right)dx}+\frac{4}{\mu}PV(t)\right)+\gamma a\int_{0}^{1}{\frac{\xi_{N,x}^{2}}{\xi^{\gamma+1}_{N}}dx}+2\int_{0}^{1}{v_{N,x}^{2}dx}\\ \leq\frac{5}{\mu}\|f_{N}\|_{L^{2}(0,1)}^{2}+\frac{|\pi_{t}|^{2}}{4}+\int_{0}^{1}{\left(\frac{\mu}{2}\xi_{N,x}^{2}+\frac{2}{\mu}v_{N}^{2}-v_{N}\xi_{N,x}+\frac{4}{\mu}G(\xi_{N}^{-1})\right)dx}+\frac{4}{\mu}PV(t).\label{gronwall0}
\end{multline}
Denote:
\begin{equation}
\eta(t):=\int_{0}^{1}{\left(\frac{\mu}{2}\xi_{N,x}^{2}+\frac{2}{\mu}v_{N}^{2}-v_{N}\xi_{N,x}+\frac{4}{\mu}G(\xi_{N}^{-1})\right)dx}+\frac{4}{\mu}PV(t),\nonumber
\end{equation}
\begin{equation}
\chi(t):=\frac{5}{\mu}\|f_{N}\|^{2}_{L^{2}(0,1)}+\frac{|\pi_{t}|^{2}}{4},\nonumber
\end{equation}
using again Cauchy's inequality we deduce
\begin{equation}
\eta(t)\geq\int_{0}^{1}{\left(\frac{\mu}{4}\xi_{N,x}^{2}+\frac{1}{\mu}v_{N}^{2}+\frac{4}{\mu}G(\xi_{N}^{-1})\right)dx}+\frac{4}{\mu}PV(t)\geq 0.\nonumber
\end{equation}
Hence, according to our denotations, inequality (\ref{gronwall0}) reads
\begin{equation}
\eta '(t)\leq \eta(t)+\chi(t).\nonumber
\end{equation}
Using the Gronwall inequality we get the following estimates:
\begin{eqnarray*}
\sup_{t\in[0,T]}\|\xi_{N}(t)\|^{2}_{H^{1}(0,1)}&\leq&\frac{4}{\mu}\mathrm{e}^{T}\left(\eta(0)+\frac{5}{\mu}\|f_{N}\|^{2}_{L_{2}(0,T;L_{2}(0,1))}+\frac{1}{4}\|\pi_{t}\|^{2}_{L_{2}(0,T)}\right),\\
\sup_{t\in[0,T]}\|v_{N}(t)\|^{2}_{L^{2}(0,1)}&\leq& \mu\mathrm{e}^{T}\left(\eta(0)+\frac{5}{\mu}\|f_{N}\|^{2}_{L_{2}(0,T;L_{2}(0,1))}+\frac{1}{4}\|\pi_{t}\|^{2}_{L_{2}(0,T)}\right),\\
\sup_{0<t<T} V(t)&\leq&\frac{\mu}{4P}\mathrm{e}^{T}\left(\eta(0)+\frac{5}{\mu}\|f_{N}\|^{2}_{L_{2}(0,T;L_{2}(0,1))}+\frac{1}{4}\|\pi_{t}\|^{2}_{L_{2}(0,T)}\right),
\end{eqnarray*}
and 
\begin{multline}
\sup_{t\in[0,T]}\frac{4}{\mu}\int_{0}^{1}{G(\xi_{N}^{-1})(x,t)dx}=\sup_{t\in[0,T]}\frac{4a}{\mu(\gamma-1)}\int_{0}^{1}{\xi_{N}^{1-\gamma}(x,t)dx}=\\
=\sup_{t\in[0,T]}\frac{4a}{\mu(\gamma-1)}\|\xi_{N}^{\frac{1-\gamma}{2}}(t)\|^{2}_{L^{2}(0,1)}\leq\mathrm{e}^{T}\left(\eta(0)+\frac{5}{\mu}\|f_{N}\|^{2}_{L_{2}(0,T;L_{2}(0,1))}+\frac{1}{4}\|\pi_{t}\|^{2}_{L_{2}(0,T)}\right).
\end{multline}
Additionally, integrating (\ref{gronwall0}) with respect to $t$ in the interval $[0,T]$ we find 
\begin{multline*}
\|v_{N,x}\|^{2}_{L_{2}(0,T;L_{2}(0,1))}=\|\xi_{N,t}\|^{2}_{L_{2}(0,T;L_{2}(0,1))}\\ \leq\frac{T\mathrm{e}^{T}+1}{2}\left(\eta(0)+\frac{5}{\mu}\|f_{N}\|^{2}_{L_{2}(0,T;L_{2}(0,1))}+\frac{1}{4}\|\pi_{t}\|^{2}_{L_{2}(0,T)}\right),
\end{multline*}
and
\begin{multline}\nonumber
\gamma a\int_{0}^{T}\!\!\!{\int_{0}^{1}{\frac{(\xi_{N,x})^{2}}{\xi^{\gamma+1}_{N}}(x,t)dx}dt}=\frac{4\gamma a}{(1-\gamma)^{2}}\int_{0}^{T}\!\!\!{\int_{0}^{1}{\left[\left(\xi_{N}^{\frac{1-\gamma}{2}}\right)_{x}\right]^{2}(x,t)dx}dt}=\\
{}=\frac{4\gamma a}{(1-\gamma)^{2}}\|\xi_{N}^{\frac{1-\gamma}{2}}\|^{2}_{L_{2}(0,T;H^{1}(0,1))}\leq(T\mathrm{e}^{T}+1)\left(\eta(0)+\frac{5}{\mu}\|f_{N}\|^{2}_{L_{2}(0,T;L_{2}(0,1))}+\frac{1}{4}\|\pi_{t}\|^{2}_{L_{2}(0,T)}\right).
\end{multline}
\begin{remark}
From the first energy estimate $f_{N}$ is bouded in $L_{\infty}(0,T;L_{2}(0,1))$, in particular
\begin{equation}\nonumber
\|f_{N}\|^{2}_{L_{2}(0,T;L_{2}(0,1))}<C_{1},
\end{equation}
where $C_{1}$ is a constant, that depends on $P,\ a,\ \mu,\  \gamma$ and initial data, but it does not depend on $T$.\\
\end{remark}
\begin{remark}
Recalling Section \ref{rozdz.3} there exists a constant $C_{2}$ depending on $P,\ a,\  \gamma$ and $\xi_{0}$, such that  $\sup_{t\in[0,\infty)}\pi_{t}(t)\leq C_{2}$. 
\end{remark}
\noindent These remarks complete the proof.
{\begin{flushright}
$\Box$
\end{flushright}}
\noindent Till now we have proved the following inclusions:
\begin{eqnarray}
\xi_{N}&\in& L_{\infty}(0,T;H^{1}(0,1)),\label{xione}\\
\xi_{N}^{\frac{1-\gamma}{2}}&\in& L_{\infty}(0,T;L_{2}(0,1)),\label{oddzielenie}\\
\xi_{N}^{\frac{1-\gamma}{2}}&\in& L_{2}(0,T;H^{1}(0,1)),\label{cisnienie}\\
\xi_{N,t}&\in& L_{2}(0,T;L_{2}(0,1)),\label{xitwo}\\
v_{N}&\in& L_{\infty}(0,T;L_{2}(0,1)),\nonumber\\
v_{N}&\in& L_{2}(0,T;H^{1}(0,1)).\nonumber
\end{eqnarray}
At the beginning of this subsection we substantiated that for $N$ sufficiently large, $\xi_{N}(x,t)>0$, now we will prove that it is indeed separated from zero.
\begin{lemat}\label{odciecie}
There exists a  positive constant K, such that for each integer $N\geq 0$
\begin{equation}\nonumber
\xi_{N}(x,t)\geq K>0\quad\  for\  (x,t)\in[0,1]\times[0,T].
\end{equation}
\end{lemat}
\dow\\
Let $\alpha>0$, then
\begin{equation}
\sup_{t\in[0,T]}\int_{0}^{1}{\left|\left(\xi_{N}^{-\alpha}\right)_{x}\right|}=
\sup_{t\in[0,T]}\int_{0}^{1}{\left|\xi_{N}\right|^{-(\alpha+1)}\left|\xi_{N,x}\right|},
\end{equation}
employing Cauchy-Schwarz inequality we obtain that for any $t\in[0,T]$
\begin{equation}\nonumber
\int_{0}^{1}{\left|\xi_{N}\right|^{-(\alpha+1)} \left|\xi_{N,x}\right|}\leq\sup_{0<t<T}
\|\xi_{N,x}\|_{L^{2}(0,1)}
\left(\int_{0}^{1}{\left|\xi_{N}\right|}^{-2(\alpha+1)}\right)^{\frac{1}{2}}.
\end{equation}
From the inclusion (\ref{oddzielenie}) we deduce
\begin{eqnarray}\nonumber
\xi_{N}^{-(\alpha+1)}\in L_{\infty}(0,T;L_{2}(0,1))
\end{eqnarray} 
iff $-2(\alpha+1)=1-\gamma$, and $\gamma>3$, then
\begin{equation}\nonumber
\xi_{N}^{-\alpha}\in L_{\infty}(0,T;W^{1}_{1}(0,1))\subset L_{\infty}(0,T;L_{\infty}(0,1)).
\end{equation}
Since
\begin{equation}
\|\xi_{N}^{-\alpha}\|_{L_{\infty}(0,T;L_{\infty}(0,1))}\leq M,
\end{equation}
for some constant $M$ and $\alpha>0$, there exists a constant $K=M^{\frac{-1}{\alpha}}$ such that
\begin{equation}\nonumber
\|\xi_{N}\|_{L_{\infty}(0,T;L_{\infty}(0,1))}\geq K>0,
\end{equation}
but we already know that  $\xi_{N}(x,t)>0$, thus
\begin{equation}\nonumber
\xi_{N}(x,t)\geq K>0.
\end{equation}
{\begin{flushright}
$\Box$
\end{flushright}}
\noindent Let estimate the pressure norm now.\\
Recall $p(\xi_{N}^{-1})=a\xi_{N}^{-\gamma}$, initially we assumed that $\gamma>1$, but from this moment we require $\gamma$ to be bigger than 3 in order to apply above lemma.\\
For such $\gamma$ the following sequence of inequalities holds
\begin{multline}\nonumber
\|p(\xi_{N}^{-1})\|_{L_{2}(0,T;H^{1}(0,1))}\leq C_{P}a\|\left(\xi_{N}^{-\gamma}\right)_{x}\|_{L_{2}(0,T;L_{2}(0,1))}\\ \leq
C_{P}aK^{\frac{-{\gamma+1}}{2}}\|\left(\xi_{N}^{\frac{1-\gamma}{2}}\right)_{x}\|_{L_{2}(0,T;L_{2}(0,1))}\leq C_{P}aK^{\frac{-{\gamma+1}}{2}}\|\xi_{N}^{\frac{1-\gamma}{2}}\|_{L^{2}(0,T;H^{1}(0,1))},
\end{multline}
where $C_{P}$ is a constant from Poincare's inequality, thus accoring to (\ref{cisnienie}) we have
\begin{equation}
p(\xi_{N}^{-1})\in L_{2}(0,T;H^{1}(0,1)).\label{normacisnienia}
\end{equation}

\subsection{Existence}
The estimates from the previous subsection imply
\begin{multline}\nonumber
\sup_{t\in[0,T]}\|v_{N}(t)\|_{L_{2}(0,1)}+\sup_{t\in[0,T]}\|\xi_{N}(t)\|_{H^{1}(0,1)}+\|v_{N}\|_{L_{2}(0,T;H^{1}(0,1))}+\|p(\xi_{N}^{-1})\|_{L_{2}(0,T;H^{1}(0,1))}\leq C.
\end{multline}
for some constant $C$ depending on $\mu$, $P$, $\gamma$, $a$, initial data and $T$. As a result we may estimate $\|v_{N,t}\|_{L_{2}(0,T;H^{-1}(0,1))}$.\\
\\
Now we will pass to limits as $N\rightarrow\infty$ to obtain a weak solution to our initial-value problem in the sense of the Definition.\\

\begin{twr}
There exists a weak solution of (\ref{eq:system}).
\end{twr}
\dow\\
Since the sequence $\{v_{N}\}_{N=1}^{\infty}$ is bounded in $L_{2}(0,T;H^{1}(0,1))$, and $\{v_{N,t}\}_{N=1}^{\infty}$ is bounded in $L_{2}(0,T;H^{-1}(0,1))$, 
 there exists a subsequence $\{v_{N_{l}}\}_{l=1}^{\infty}\subset\{v_{N}\}_{N=1}^{\infty}$, such that
\begin{displaymath}
\begin{array}{rcrl}
v_{N_{l}} & \rightharpoonup&  v\in& L_{2}(0,T;H^{1}(0,1)),\\
v_{N_{l},t} & \rightharpoonup&  v_{t}\in& L_{2}(0,T;H^{-1}(0,1)).
\end{array}
\end{displaymath} 
With the same manner we can conclude that 
\begin{eqnarray*}
\xi_{N_{l_{k}}}&\rightharpoonup^{*}&\xi\in L_{\infty}(0,T;H^{1}(0,1)),\\
\xi(x,t)&=&\pi(t)\quad\mathrm{on}\quad \{0\}\times [0,T]\quad\mathrm{and} \quad \{1\}\times[0,T],
\end{eqnarray*}
for some subsequence
$\left\{\xi_{N_{l_{k}}}\right\}_{k=1}^{\infty}\subset\left\{\xi_{N_{l}}\right\}_{l=1}^{\infty}$, 
starting from here we will be calling this subsequence $\xi_{N}$, $v_{N}$.\\
According to the inclusions (\ref{xione}), (\ref{xitwo}) 
\begin{equation*}
\xi_{N}\in H^{1}((0,1)\times(0,T)),
\end{equation*}
hence by the Rellich-Kondrachov Compactness Theorem we obtaine strong convergence of some subsequence $\left\{\xi_{N_{l}}\right\}_{l=1}^{\infty}\subset\left\{\xi_{N}\right\}_{N=1}^{\infty}$, $\xi_{N_{l}}\rightarrow \xi$ in $L_{2}((0,1)\times (0,T))$.\\
Next, observe that since $p(\xi_{N}^{-1})$ is bouded in ${L_{2}(0,T;H^{1}(0,1))}$ there exists a weakly convergent subsequence $\left\{\xi_{N_{l_{k}}}\right\}_{k=1}^{\infty}\subset\left\{\xi_{N_{l}}\right\}_{l=1}^{\infty}$ to some function in $L_{2}(0,T;H^{1}(0,1))$, i.e.
\begin{equation}
p\left(\xi_{N_{l_{k}}}^{-1}\right)\rightharpoonup\overline{p\left(\xi_{N_{l_{k}}}^{-1}\right)}\in L_{2}(0,T;H^{1}(0,1)).\nonumber
\end{equation}
\begin{lemat}\label{jegorow}
Providing $p(\xi^{-1})$ is a continuous function of $\xi$ and that with an accuracy to subsequence
\begin{eqnarray*}
\xi_{N_{l}}&\rightarrow&\xi\quad strongly\  in\   L_{2}((0,1)\times (0,T)),\\
p\left(\xi_{N_{l_{k}}}^{-1}\right)&\rightharpoonup& \overline{p\left(\xi_{N_{l_{k}}}^{-1}\right)}\quad weakly\  in\  L_{2}(0,T;H^{1}(0,1)),
\end{eqnarray*}
then $\overline{p(\xi_{N_{l_{k}}}^{-1})}=p(\xi^{-1})$ holds a.e. in $(0,1)\times[0,T]$.
\end{lemat}
\noindent The Proof of this Lemma follows easily from Egoroff's Theorem.\\

\noindent Similaryly, since  $\left\{\xi_{N}^{\frac{1-\gamma}{2}}\right\}_{N=1}^{\infty}$ is bounded in
$L_{\infty}(0,T;L_{2}(0,1))$, we may choose such subsequence $\left\{\xi_{N_{l}}^{\frac{1-\gamma}{2}}\right\}_{l=1}^{\infty}\subset\left\{\xi_{N}^{\frac{1-\gamma}{2}}\right\}_{N=1}^{\infty}$ that
\begin{equation}
\xi_{N_{l}}^{\frac{1-\gamma}{2}}\rightharpoonup^{*}\overline{\xi^{\frac{1-\gamma}{2}}_{N_{l}}}\in L_{\infty}(0,T;L_{2}(0,1))\nonumber
\end{equation}
and by the same argument prove
\begin{equation}
\overline{\xi^{\frac{1-\gamma}{2}}_{N_{l}}}=\xi^{\frac{1-\gamma}{2}}\nonumber
\end{equation}
a.e. in $(0,1)\times[0,T]$.\\
Repeating the procedure from proof of Lemma (\ref{odciecie}) for $\xi(x,t)$ we may show
\begin{equation}
\xi(x,t)\geq C>0\nonumber
\end{equation}
for  $(x,t)\in(0,1)\times[0,T].$\\
\\
According to our previous remarks it is possible to pass to limits in the weak formulation, and by the density argument we get that
\begin{equation}
(v_{t},\varphi)-(p(\xi^{-1})-P,\varphi_{x})+\mu(v_{x},\varphi_{x})+(f,\varphi)=0\nonumber
\end{equation}
holds for each $\varphi\in W$ in the sense of distributions on time interval $[0,T]$, moreover
\begin{equation}
v_{x}(x,t)=\xi_{t}(x,t)\quad \mathrm{on\ } (0,1)\times[0,T].\nonumber
\end{equation}

\noindent In order to prove $v(x,0)=v_{0}$, it suffices to take as a test function $\varphi\in W$ such that $\varphi(x,T)=0$  and pass to the weak limits taking into account that $v_{N_{l}}(0)\rightarrow v_{0}$ in $L_{2}(0,1)$.\\
Equation $v_{x}(x,t)=\xi_{t}(x,t)\  \mathrm{on\ } (0,1)\times[0,T]$ enables to prove a suitable initial condition for $\xi(x,t)$
\begin{equation}
\xi(x,0)=\xi_{0}. \nonumber
\end{equation}
To complete the prove of existence of weak solutions there is a need to show higher regularity of $v$.
\begin{lemat}
\begin{equation}\nonumber
v\in W^{2,1}_{2}((0,1)\times(0,T)).
\end{equation}
\end{lemat}
\dow\\
Multiply equation (\ref{gamma}) by $\alpha^{N}_{k,t}$ and sum $k=1,\ldots,N$, to discover 
\begin{equation}\nonumber
(v_{N,t},\tilde{v}_{N,t})-(p(\xi_{N}^{-1})-P,\tilde{v}_{N,tx})+\mu(\mathbf{T}v_{N,x},\tilde{v}_{N,tx})=0,
\end{equation}
where $\tilde{v}_{N}(x,t)=\sum_{k=1}^{N}\alpha^{N}_{k}(t)w_{k}(x)$.\\
Let $u(x,t)=-\frac{\pi_{t}(t)}{2}+x\pi_{t}(t)$, $u_{t}(x,t)\in W^{N}$, thus it satisfies
\begin{equation}\nonumber
(v_{N,t},u_{t})-(p(\xi_{N}^{-1})-P,u_{tx})+\mu(\mathbf{T}v_{N,x},u_{tx})=0.
\end{equation}
Adding these two equalities and using the formula (\ref{laplasjan}) we get
\begin{equation}
\int_{0}^{1}(v_{N,t})^{2}dx-\mu\int_{0}^{1}\tilde{v}_{N,xx}(\tilde{v}_{t}+u_{t})dx=-\int_{0}^{t}p\left(\xi_{N}^{-1}\right)_{x}v_{N,t}
-\mu\int_{0}^{1}(1-\mathbf{T})v_{N,xx}v_{N,t},\label{regularnoscv}
\end{equation}
since $v_{N}(x,t)=u(x,t)+\tilde{v}(x,t)$, in particular $v_{N,xx}=\tilde{v}_{N,xx}$.
Now employ Caychy's inequality and integrate by parts to find
\begin{equation}\nonumber
\frac{\|v_{N,t}\|^{2}_{L_{2}(0,1)}}{3}+\frac{\mu}{2}\frac{d}{dt}\int_{0}^{1}(\tilde{v}_{N,x})^{2}dx+ \mu\int_{0}^{1}\tilde{v}_{N,x}u_{tx}\leq
\frac{3}{4}\|p\left(\xi_{N}^{-1}\right)_{x}\|_{L_{2}(0,1)}^{2}+\frac{3\mu^{2}}{4}\|(1-\mathbf{T})v_{N,xx}\|^{2}_{L_{2}(0,1)}
.\end{equation}
Note $\|\tilde{v}_{N,x}\|_{L_{2}(0,T;L_{2}(0,1))}\leq\|v_{N}\|_{L_{2}(0,T;H^{1}(0,1))}$, and $u_{tx}=\pi_{tt}\in L_{2}(0,T;L_{2}(0,1))$, therefore integreting with respect to $t$ in the interval  $[0,T]$ we conclude
\begin{equation}\nonumber
v_{N,t}\in L_{2}(0,T;L_{2}(0,1)),
\end{equation}
also by (\ref{regularnoscv}) and Cauchy-Schwarz inequality
\begin{equation}\nonumber
v_{N,xx}\in L_{2}(0,T;L_{2}(0,1)).
\end{equation}
As $N$ approaches infinity we obtain required smoothness of  $v$, thus the proof of existence is complete. 
\begin{flushright}
$\Box$
\end{flushright}

\subsection{Uniqueness}
\begin{twr}
A weak solution of (\ref{eq:system}) is unique.
\end{twr}
\dow\\
Let assume that there are two weak solutions $(v_{1},\xi_{1})$ and $(v_{2},\xi_{2})$ of the system (\ref{eq:system}) in the sense of the Definition.\\
Denote $\omega=v_{1}-v_{2}$, $\psi=\xi_{1}-\xi_{2}$ and insert $\omega\in W$ in place of the test function $\varphi$ in the definition of weak solution
\begin{equation}\nonumber
(\omega_{t},\omega)+(\left( p(\xi_{1}^{-1})-p(\xi_{2}^{-1})\right)_{x},\omega)=(\mu(\mathbf{T}\omega)_{xx},\omega).
\end{equation}
Employing the formula
\begin{equation}\nonumber
f(\xi_{1})-f(\xi_{2})=(\xi_{1}-\xi_{2})\int_{0}^{1}f'(s\xi_{1}+(1-s)\xi_{2})ds,
\end{equation}
for $f(\xi)=p(\xi^{-1})=\frac{a}{\xi^{\gamma}}$, integrating by parts and replacing $\omega_{x}$ by $\psi_{t}$, we find
\begin{displaymath}
\frac{d}{dt}\int_{0}^{1}\left(\frac{\omega^{2}}{2}-\frac{\psi^{2}}{2}\int_{0}^{1}f'(s\xi_{1}+(1-s)\xi_{2})ds\right)dx=~~~~~~~~~~~~~~~~~~~~~~~~~~~~~~~~
\end{displaymath}
\vskip-18pt
\begin{eqnarray*}
~~~~~~~~~~~~~~~~~~~~~~~~~~~~~~~~~~~~~~~~~~~~~~~~~~~~{}&=&-\int_{0}^{1}\frac{\psi^{2}}{2}\int_{0}^{1}f''(s\xi_{1}+(1-s)\xi_{2})ds\xi_{t}dx-\mu\int_{0}^{1}(\mathbf{T}\omega_{x})^{2}dx\\
{}&\leq& C(t)\|\xi_{t}(t)\|_{L_{\infty}(0,1)}\int_{0}^{1}\frac{\psi^{2}}{2}dx
,\end{eqnarray*}
where $C(t)=\|\int_{0}^{1}f''(s\xi_{1}+(1-s)\xi_{2})ds(t)\|_{L_{\infty}(0,1)}$, $\sup_{t\in[0,T]}C(t)<\infty$, and
in accordance with  $\xi_{t}=v_{x}\in L_{2}(0,T;H^{1}(0,1))\subset L_{1}(0,T;L_{\infty}(0,1))$.\\
Since $f(\xi)=\frac{a}{\xi^{\gamma}}$ is monotonically increasing function, its derivative is strictly negative, thus
\begin{equation*}
\frac{\omega^{2}}{2}-\frac{\psi^{2}}{2}\int_{0}^{1}f'(s\xi_{1}+(1-s)\xi_{2})ds\geq 0.
\end{equation*}
Let $\phi(t)=C(t)\|\xi_{t}(t)\|_{L_{\infty}(0,1)}$, note $\phi\in L_{1}(0,T)$, then by Gronwall's inequality and the initial conditions $\omega(x,0)=0,\ \psi(x,0)=0$ we discover
\begin{equation}\nonumber
\int_{0}^{1}\left(\frac{\omega^{2}}{2}-\frac{\psi^{2}}{2}\int_{0}^{1}f'(s\xi_{1}+(1-s)\xi_{2})ds\right)dx\leq 0.
\end{equation}
and therefore $\psi(x,t)=\omega(x,t)\equiv 0$. 
\begin{flushright}
$\Box$
\end{flushright}

\section{Proof of Theorem B}
To obtain a global in time existence in case when  the local existence has been already proved we have to show only some uniform in time estimates for solutions of (\ref{eq:system}).
\begin{lemat}
For a solution of (\ref{eq:system}) we have
\begin{multline*}
\int_{0}^{1}{\left(\frac{1}{2}v^{2}(x,t)+G(\xi^{-1})(x,t)\right)dx}+PV(t)+\mu\int^{t}_{0}\!\!\!{\int_{0}^{1}{|\mathbf{T}v_{x}|^{2}dx}ds}=\\
=\int_{0}^{1}{\left(\frac{1}{2}v_{0}^{2}+G(\xi_{0}^{-1})\right)dx}+PR(0).
\end{multline*}
\end{lemat}
\dow\\
Multiplying the first equation of (\ref{eq:system}) by $v$, integrating over $[0,1]$ and repeating the proof of Lemma \ref{Ienerget} we get
\begin{equation}
\frac{d}{dt}\int_{(0,1)}{\left(\frac{1}{2}v^{2}+G(\xi^{-1})\right)dx}+\frac{d}{dt}PV(t)+\mu\int_{0}^{1}{|\mathbf{T}v_{x}|^{2}dx}=0,\label{globalne1}
\end{equation}
integrating over $[0,t]$ we complete the proof.
\begin{flushright}
$\Box$
\end{flushright}
Therefore we have
\begin{eqnarray}
v&\in& L_{\infty}(0,\infty;L_{2}(0,1)),\label{I}\\
G(\xi^{-1})&\in& L_{\infty}(0,\infty;L_{1}(0,1)),\label{II}\\
V(t)&\in& L_{\infty}(0,\infty),\label{III}\\
\mathbf{T}v_{x}&\in& L_{2}(0,\infty;L_{2}(0,1)).\nonumber
\end{eqnarray}
Now we will show that for $\xi$ holds 
\begin{equation}\nonumber
\xi\in L_{\infty}(0,\infty;H^{1}(0,1)).
\end{equation}

\begin{lemat}[The upper bound on $\xi$]
If $\frac{P}{\mu}$ is sufficiently large, then
\begin{equation}\nonumber
\xi(x,t)\leq\max\left\{\xi_{min},\sup_{x\in(0,1)}(\xi_{0}+\frac{1}{\mu}U(x,0))+\frac{1}{\mu}M_{U}-t\frac{P}{4\mu}    \right\}
\end{equation}
for $t\in [0,T_{max})$, where $T_{max}$ is the maximal time of existence of solutions in sense of the Definition.
\end{lemat}
\dow\\
\noindent Introduce a new function $U(x,t)$ as follows
\begin{equation}\nonumber
U(x,t)=\int_{0}^{x}{v(s,t)ds}.
\end{equation}
Since $\int_{0}^{1}{v(x,t)dx}=0$, we have
\begin{equation}\nonumber
U(0,t)=0,\quad U(1,t)=0,
\end{equation}
and
\begin{equation}\nonumber
(v_{t},\varphi)=(U_{xt},\varphi)=-(U_{t},\varphi_{x}).
\end{equation}
Therefore we rewrite the weak formulation in the form
\begin{equation}
\left(\left(\xi-\frac{1}{\mu}U\right)_{t},\varphi_{x}\right)=
\left(\frac{1}{\mu}(p(\xi^{-1})-P)+(1-\mathbf{T})\xi_{t},\varphi_{x}\right).\label{zapis2}
\end{equation}
Moreover, by (\ref{I})
\begin{equation}\nonumber
U\in L_{\infty}(0,\infty;H^{1}_{0}(0,1))\subset L_{\infty}((0,\infty)\times(0,1)),
\end{equation}
and by properties of operator $\mathbf{T}$
\begin{equation}\nonumber
(1-\mathbf{T})v_{x}\in L_{\infty}((0,\infty)\times(0,1)).
\end{equation}
In particular, the following bounds are true
\begin{equation}\nonumber
\|U\|_{L_{\infty}((0,\infty)\times(0,1))}\leq M_{U},\qquad\|(1-T)\xi_{t}\|_{L_{\infty}((0,\infty)\times(0,1))}\leq C.
\end{equation}
The constant C is independent of $\mu$, thus we choose $\mu$ sufficiently small to keep
\begin{equation}\nonumber
\frac{-P}{\mu}+(1-T)\xi_{t}\leq \frac{-P}{2\mu}.
\end{equation}
Let $\xi_{min}$ be a positive constant which satisfies 
\begin{equation}\nonumber
p(\xi_{min}^{-1})<\frac{P}{4},
\end{equation}
by the Lemma \ref{nabrzegu} it can be done for any case. Then we see that 
\begin{equation}\nonumber
\frac{1}{\mu}p(\xi_{min}^{-1})-\frac{P}{\mu}+(1-T)\xi_{t}\leq \frac{-P}{4\mu}.
\end{equation}
Let
\begin{equation}\nonumber
N=\left\{(x,t):\xi(x,t)-\frac{1}{\mu}U(x,t)\geq\xi_{min}-\frac{1}{\mu}M_{U}\right\},
\end{equation}
then remembering that $p(\cdot)$ is an increasing function, we get $\xi|_{N}\geq\xi_{min}$.\\
Since we require only that  $\varphi\in W$, i.e. $\int_{0}^{1}{\varphi(x,t)}=0$, 
hence there is no restriction on $\varphi_{x}$ which appears in (\ref{zapis2}). In particular, taking $\varphi$ such that $\mathrm{supp}\  \varphi_{x}\subset N$, we conclude
\begin{equation}
\frac{\partial}{\partial t}\left(\xi-\frac{1}{\mu}U\right)\Big|_{N}\leq\frac{-P}{4\mu}.\label{Tmin}
\end{equation}
Then there are two possibilities:
\begin{enumerate}
\item If 
\begin{equation}\nonumber
\sup_{x\in(0,1)}\xi_{0}(x)-\frac{1}{\mu}U(x,0)\leq\xi_{min}-\frac{1}{\mu}M_{U},
\end{equation}
then if there was
\begin{equation}
\sup_{(x,t)\in(0,1)\times(0,\infty)}\xi(x,t)-\frac{1}{\mu}U(x,t)>\xi_{min}-\frac{1}{\mu}M_{U}\label{sprzecznosc}
\end{equation}
 it would exist a point $(x_{0}, t_{0})$ such that
\begin{equation}\nonumber
\xi(x_{0},t_{0})-\frac{1}{\mu}U(x_{0},t_{0})=\xi_{min}-\frac{1}{\mu}M_{U}
\end{equation}
and by (\ref{Tmin})
\begin{equation}\nonumber
\frac{\partial}{\partial t}\left(\xi(x_{0},t_{0})-\frac{1}{\mu}U(x_{0},t_{0}) \right)\leq\frac{-P}{4\mu},
\end{equation}
but it is a contradition to (\ref{sprzecznosc}).\\
Thus in this case the following bound is valid
\begin{equation}\nonumber
\sup_{(x,t)\in(0,1)\times(0,\infty)}\xi(x,t)-\frac{1}{\mu}U(x,t)\leq\xi_{min}-\frac{1}{\mu}M_{U}.
\end{equation}
\item If
\begin{equation}\nonumber
\sup_{x\in(0,1)}\xi_{0}(x)-\frac{1}{\mu}U(x,0)>\xi_{min}-\frac{1}{\mu}M_{U},
\end{equation}
then for $t\in[0,T_{min})$ we have (\ref{Tmin}), it means that the function $\xi(x,t)-\frac{1}{\mu}U(x,t) $ decreases untill it reaches the value $\xi_{min}-\frac{1}{\mu}M_{U}$. We compute $T_{min}$ from the condition (\ref{Tmin})
\begin{equation}\nonumber
T_{min}=\frac{4\mu}{P}\left( \sup_{x\in(0,1)}(\xi_{0}-\frac{1}{\mu}U(x,0))-\xi_{min}+\frac{1}{\mu}M_{U}\right).
\end{equation} 
Thus for $t\in[0,T_{min}]$ we have
\begin{equation}\nonumber
\xi(x,t)\leq\sup_{x\in(0,1)}(\xi_{0}-\frac{1}{\mu}U(x,0))+\frac{1}{\mu}M_{U}-t\frac{P}{4\mu}.
\end{equation}
and for $t\in(T_{min},\infty)$ the bound from previous case is valid 
\end{enumerate}
Combining these two case we complete the proof.
\begin{flushright}
$\Box$
\end{flushright}
So, there exists a constant $\xi_{+}$ such that
\begin{equation*}
\xi(x,t)\leq\xi_{+}
\end{equation*}
for all $(x,t)\in(0,1)\times[0,\infty)$.
\begin{lemat}
For $\mu\leq\frac{a\gamma}{\xi_{+}
^{\gamma+1}}$ there exists a constant $K$ depending on $P, a,\ \gamma$ and the initial data, such that
\begin{equation*}
\sup_{t\in[0,\infty)}\|\xi\|_{H^{1}(0,1)}\leq K.
\end{equation*}
\end{lemat}
\dow\\
\noindent 
Multiplying the first equation of (\ref{eq:system}) by $\xi_{x}$, integrating over $x\in[0,1]$ and  rearranging it the same way we did proving the Lemma (\ref{IIenerget}), we find
\begin{multline}\label{globalne2}
\frac{d}{dt}\int_{0}^{1}{\left(\frac{\mu}{2}\xi_{x}^2-v\xi_{x}\right)dx}+\gamma a\int_{0}^{1}{\frac{\xi_{x}^{2}}{\xi^{\gamma+1}}dx}=\int_{0}^{1}{v_{x}^{2}dx}+\int_{0}^{1}{f\xi_{x}dx}-\pi_{t}\int_{0}^{1}v_{x}dx,
\end{multline}
where $f(x,t)=\mu(1-\mathbf{T})v_{xx}\in L_{\infty}(0,\infty;L_{2}(0,1))$.\\
Now we multiply equation (\ref{globalne1})  by  $\frac{4}{\mu}$ and add to (\ref{globalne2}), then recalling the formula $
\mu\int_{0}^{1}|Tv_{x}|^{2}dx=\mu\int_{0}^{1}|v_{x}|^{2}dx+\int_{0}^{1}fvdx$, we get
 \begin{multline}\label{globalne3}
\frac{d}{dt}\left(\int_{0}^{1}{\left(\frac{\mu}{2}\xi_{x}^{2}+\frac{2}{\mu}v^{2}-v\xi_{x}+\frac{4}{\mu}G(\xi^{-1})\right)dx}+
\frac{4}{\mu}PV(t)\right)+\gamma a\int_{0}^{1}{\frac{\xi_{x}^{2}}{\xi^{\gamma+1}}dx}+3\int_{0}^{1}{v_{x}^{2}dx}=\\
=\int_{0}^{1}{(f\xi_{x}-\frac{4}{\mu}fv)dx}-\pi_{t}\int_{0}^{1}v_{x}dx.
\end{multline}
It follows from the previous lemma that
\begin{equation}
\int_{0}^{1}{\frac{(\xi_{x})^{2}}{\xi^{\gamma+1}}dx}\geq \frac{1}{\xi_{+}^{\gamma+1}}\int_{0}^{1}\xi_{x}^{2}.\nonumber
\end{equation}
Since $\mu<\frac{\gamma a}{\xi_{+}^{\gamma+1}}$, by H\"{o}lder's, Cauchy's and Poincare's (with a constant $C_{P}$) inequalities we obtain from (\ref{globalne3}) an expression
\begin{multline}\label{globalne5}
\frac{d}{dt}\left(\int_{0}^{1}{\left(\frac{\mu}{2}\xi_{x}^2+\frac{2}{\mu}v^{2}-v\xi_{x}+\frac{4}{\mu}G(\xi^{-1})\right)dx} +\frac{4}{\mu}PV(t)\right)+\int_{0}^{1}{\left(\frac{\mu}{2}\xi_{x}^2+\frac{2}{\mu}v^{2}-v\xi_{x}\right)dx}\\
\leq\left(\frac{1}{\mu}+1\right)\|f\|_{L_{2}(0,1)}^{2}+\frac{1}{4}|\pi_{t}|^{2}+\left(\frac{2}{C_{P}}-\frac{4}{\mu}\right)\int_{0}^{1}v^{2}dx.
\end{multline}
Denote:
\begin{equation*}
\eta(t)=\int_{0}^{1}{\left(\frac{\mu}{2}\xi_{x}^2+\frac{2}{\mu}v^{2}-v\xi_{x}+\frac{4}{\mu}G(\xi^{-1})\right)dx} +\frac{4}{\mu}PV(t),
\end{equation*}
then (\ref{globalne5}) reads
\begin{equation*}
\frac{d}{dt}\eta(t)+\eta(t)\leq\left(\frac{1}{\mu}+1\right)\|f\|_{L_{2}(0,1)}^{2}+\frac{1}{4}|\pi_{t}|^{2}+\left(\frac{2}{C_{P}}-\frac{4}{\mu}\right)\int_{0}^{1}v^{2}dx
+\int_{0}^{1}\frac{4}{\mu}G(\xi^{-1})dx+\frac{4}{\mu}PV(t).
\end{equation*}
The right hand side of this inequality is bounded by a constant $M$ in $L_{\infty}(0,\infty)$, which is a coclusion from (\ref{I}), (\ref{II}), (\ref{III}).
In accordance with our denotations the following bound holds
\begin{equation*}
\eta(t)\leq\eta(0)\mathrm{e}^{-t}+M\mathrm{e}^{-t}+M,
\end{equation*}
thus
\begin{equation*}
\sup_{t\in[0,\infty)}\eta(t)\leq\eta(0)+2M.
\end{equation*}
Observe, taht by  Cauchy's inequality $\eta(t)\geq\int_{0}^{1}{\left(\frac{\mu}{4}\xi_{x}^2+\frac{1}{\mu}v^{2}+\frac{4}{\mu}G(\xi^{-1})\right)dx} +\frac{4}{\mu}PV(t)$, therefore we truly have
\begin{equation*}
\xi\in L_{\infty}(0,\infty;H^{1}(0,1)).\ 
\end{equation*}
\begin{flushright}
$\Box$
\end{flushright}
The inclusion (\ref{II}) implies
\begin{equation*}
\xi^{\frac{1-\gamma}{2}}\in L_{\infty}(0,\infty;L_{2}(0,1)),
\end{equation*} 
since we already know that $\xi\in L_{\infty}(0,\infty;H^{1}(0,1))$  we may repeat an argument from Lemma \ref{odciecie} to obtain the existence of a constant $\xi_{-}>0$ such that
\begin{equation}\nonumber
\xi(x,t)\geq\xi_{-}>0
\end{equation}
for all $(x,t)\in(0,1)\times[0,\infty)$.
With this fact the proof of global in time existence is complete.

\end{document}